\documentclass[11pt]{amsart}
\usepackage{amscd,amssymb,amsmath}
\newtheorem{theorem}{Theorem}[section]

\theoremstyle{definition}
\newtheorem{definition}[theorem]{Definition}
\newtheorem{example}[theorem]{Example}
\newtheorem{proposition}[theorem]{Proposition}

\date \kill
\theoremstyle{remark}
\newtheorem{remark}{\bf Remark}[section]

\numberwithin{equation}{section}
\begin{document}
\title{Equivariant evaluation subgroups and Rhodes groups}
\author{Marek Golasi\'nski}
\address{Faculty of Mathematics and Computer Science, Nicolaus Copernicus University, Chopina 12/18, 87-100, Toru\'n, Poland}
\email{marek@mat.uni.torun.pl}
\author{Daciberg Gon\c calves}
\address{Dept. de Matem\'atica - IME - USP, Caixa Postal 66.281 - CEP 05311-970,
S\~ao Paulo - SP, Brasil}
\email{dlgoncal@ime.usp.br}
\author{Peter Wong}
\address{Department of Mathematics, Bates College, Lewiston,
ME 04240, U.S.A.}
\email{pwong@bates.edu}
\thanks{The authors are grateful to the referee for carefully reading an earlier
version of the paper and all his suggestions that make the paper more clear and readable.
This work was conducted during the second and third authors' visit to the Faculty
of Mathematics Computer Science, Nicolaus Copernicus University August 4 - 13, 2005.
The second and third authors would like to thank the Faculty of Mathematics and
Computer Science for its hospitality and support}

\begin{abstract}
In this paper, we define equivariant evaluation subgroups of the higher
Rhodes groups and study their relations with Gottlieb-Fox groups. 
\end{abstract}
\date{}
\maketitle

\newcommand{\af}{\alpha}
\newcommand{\et}{\eta}
\newcommand{\ga}{\gamma}
\newcommand{\ta}{\tau}
\newcommand{\ph}{\varphi}
\newcommand{\bt}{\beta}
\newcommand{\lb}{\lambda}
\newcommand{\wh}{\widehat}
\newcommand{\wt}{\widetilde}
\newcommand{\sg}{\sigma}
\newcommand{\om}{\omega}
\newcommand{\cH}{\mathcal H}
\newcommand{\cF}{\mathcal F}
\newcommand{\N}{\mathcal N}
\newcommand{\R}{\mathcal R}
\newcommand{\Ga}{\Gamma}
\newcommand{\cc}{\mathcal C}

\newcommand{\bea} {\begin{eqnarray*}}
\newcommand{\beq} {\begin{equation}}
\newcommand{\bey} {\begin{eqnarray}}
\newcommand{\eea} {\end{eqnarray*}}
\newcommand{\eeq} {\end{equation}}
\newcommand{\eey} {\end{eqnarray}}

\newcommand{\ovl}{\overline}
\newcommand{\vv}{\vspace{4mm}}
\newcommand{\lra}{\longrightarrow}

\bibliography{ref}
\bibliographystyle{amsplain}

\section*{Introduction}

While the Gottlieb groups $G_n(X)$ of a $CW$-space $X$ with a basepoint are important objects of study in homotopy theory,
the first Gottlieb group $G_1(X)$ actually originated from Nielsen fixed point theory. In fact, $G_1(X)$ was
first introduced by B. Jiang and it is also known as the Jiang subgroup $J(X)$.

In \cite{wong}, the classical Nielsen fixed point theory was generalized to the equivariant setting under
the presence of a group action. Subsequent works related to equivariant fixed point theory include
\cite{fagundes-goncalves} and \cite{guo}. To facilitate computation of certain equivariant Nielsen
type numbers, an equivariant Jiang condition was introduced in \cite{wong}. This condition was studied and
relaxed by Fagundes and Gon\c calves \cite{fagundes-goncalves} who gave an example in which all
equivariant fixed point classes have the same index while the space does not satisfy the equivariant Jiang
conidtion. In \cite{gg} equivariant Gottlieb groups, which are analogs of the equivariant Jiang subgroups
in higher homotopy groups, were defined and used to compute the Gottlieb groups of orbit
spaces.

This paper is organized as follows. In Section 1, we recall the definition
of the higher Rhodes groups $\sigma_n(X,x_0,G)$ and prove the following.

{\bf Theorem 1.2.} {\em Suppose a group $G$ acts freely on $X$ with a basepoint $x_0$.
Then, $\sigma_n(X,x_0,G) \stackrel{p_*}{\to} \tau_n(X/G,p(x_0))$
is an isomorphism for all $n\ge 1$.}

\noindent
This generalizes the result of Rhodes when $n=1$. Moreover, we establish a split exact sequence, similar to that proved by Fox for the torus
homotopy groups, for the Rhodes groups. We then introduce new equivariant Gottlieb groups (in contrast to
those introduced in \cite{gg}) as subgroups of the higher Rhodes groups and discuss some basic properties
in Section 2. We prove

{\bf Proposition 2.4.} {\em A space $X$ is Gottlieb if and only if it is a
Gottlieb-Fox space.}

Relationships among various evaulation subgroups are discussed in
Section 3. In particular, $n$-Gottlieb and equivariant $n$-Gottlieb spaces are related
as follows.

{\bf Theorem 3.1.} {\em Suppose a group $G$ acts freely on a spaces $X$ with
a basepoint $x_0$.
\begin{enumerate}
 \item For $n\ge 2$, if $X$ is equivariant $n$-Gottlieb then $X/G$ is $n$-Gottlieb.
 \item For $n\ge 1$, if $X/G$ is $n$-Gottlieb then $X$ is equivariant $n$-Gottlieb.
 \item Suppose $X$ is a finite aspherical $G$-CW space. If $X$ is equivariant $1$-Gottlieb then $X/G$ is $1$-Gottlieb.
\end{enumerate}}

\section{Higher Rhodes groups}

Throughout, $G$ will denote a finite group acting on a compactly generated Hausdorff path connected space $X$
with a basepoint. The associated pair $(X,G)$ is called in the literature a transformation group.

In \cite{rhodes1}, F.\ Rhodes introduced the notion of the fundamental
group $\sigma (X,x_0,G)$ of the pair $(X,G)$, where $x_0$ is a basepoint
in $X$. A typical element in $\sigma (X,x_0,G)$ is the homotopy class
$[\alpha; g]$ consisting of a path $\alpha$ in $X$ and a group element $g\in G$
such that $\alpha(0)=x_0, \alpha(1)=gx_0$. The multiplication in $\sigma (X,x_0,G)$
is given by
$$
[\alpha_1;g_1]*[\alpha_2;g_2]:=[\alpha_1+g_1\alpha_2;g_1g_2].
$$

It is easy to see that the groups $\pi_1(X,x_0)$, $\sigma (X,x_0,G)$ and $G$ fit into the following short exact sequence
$$
1\to \pi_1(X,x_0) \to \sigma (X,x_0,G) \to G \to 1.
$$
Then, for $n\ge 1$, F.\ Rhodes \cite{rhodes2} defined higher homotopy groups $\sigma_n(X,x_0,G)$ of $(X,G)$ which is
an extension of $\tau_n(X,x_0)$ by $G$ so that
\begin{equation}\label{rhodes-exact}
1\to \tau_n(X,x_0) \to\sigma_n (X,x_0,G) \to G \to 1
\end{equation}
is exact. Here, $\tau_n(X,x_0)$ denotes the $n$-th torus homotopy group of $X$ introduced
by R.\ Fox \cite{fox}. The group $\tau_n=\tau_n(X,x_0)$ is defined to be the fundamental group
of the function space $X^{T^{n-1}}$ and is uniquely determined by the groups $\tau_1, \tau_2, \ldots, \tau_{n-1}$
and the Whitehead products, where $T^{n-1}$ is the $(n-1)$-dimensional torus. The group $\tau_n$ is
non-abelian in general.

\begin{definition}
Suppose that $X$ is a $G$-space with a basepoint $x_0\in X$.
Let $C_n=I\times T^{n-1}$.  We say that a map $f : C_n\to X$ is of {\em order}
$g\in G$ provided $f(0,t_2,\ldots,t_n)=x_0$ and $f(1,t_2,\ldots,t_n)=g(x_0)$ for $(t_2,\ldots,t_n)\in T^{n-1}$. Denote by $[f;g]$
the homotopy class of a map $f: C_n\to X$ of order $g$ and by $\sigma_n(X,x_0,G)$ the set of all such homotopy classes.
We define an operation $*$ similar to the one on $\sigma (X,x_0,G)$ on the set $\sigma_n(X,x_0,G)$, i.e.,
$$[f';g']*[f;g]:=[f' + g'f;g'g].$$
This operation makes $\sigma_n(X,x_0,G)$ a group. We write $\sigma_1(X,x_0,G)$ for $\sigma(X,x_0,G)$.
Now, the orbit map $p : X\to X/G$ leads to a homomorphism $p_\ast : \sigma_n(X,x_0,G)\to \tau_n(X/G,p(x_0))$,
where $p_\ast([f;g])$ is represented by the adjoint of the composite $I\to (X/G)^{T^{n-1}}$ of $C_n=I\times T^{n-1}\stackrel{f}{\to}X\stackrel{p}{\to}X/G$
for $[f;g]\in \sigma_n(X,x_0,G)$.
\end{definition}

\bigskip

\par It was shown in \cite{rhodes1} that $p_\ast: \sigma_1(X,x_0,G)\to \tau_1(X/G, p(x_0))$
is an isomorphism for a free $G$-action on $X$. More generally, we have the following

\begin{theorem}
Suppose $G$ acts freely on $X$ with a basepoint $x_0$. Then,
$$
\sigma_n(X,x_0,G) \stackrel{p_*}{\to} \tau_n(X/G,p(x_0))
$$
is an isomorphism for all $n\ge 1$.\label{T}
\end{theorem}
\begin{proof} Since the $G$-action on $X/G$ is trivial, it follows that
$$
\sigma_n(X/G,p(x_0),G)\cong \tau_n(X/G,p(x_0)).$$
By induction, we assume that
$$
\sigma_k(X,x_0,G)\cong \tau_k(X/G,p(x_0))$$
for $k<n$. We define the projection $\sigma_n \to \sigma_{n-1}$ by $[f;g] \mapsto [\hat f;g]$,
where $\hat f:C_{n-1}\to X$ is given by $f\circ i_n$ with $i_n:C_{n-1} \hookrightarrow C_n$
defined by $i_n(t_1,t_2,\ldots,t_{n-1})=(t_1,t_2,\ldots,t_{n-1},0)$ for $(t_1,\ldots,t_{n-1})\in C_{n-1}$.
The projection $\sigma_n \to \sigma_{n-1}$ has a section $\sigma_{n-1}\to \sigma_n$ via
the projection $C_n\to C_{n-1}$ defined by $(t_1,t_2,\ldots,t_n) \mapsto (t_1,t_2,\ldots,t_{n-1})$
for $(t_1,\ldots,t_n)\in C_n$.
Consider the following commutative diagram
\begin{equation}\label{rhodes-fox-split}
\begin{CD}
1 @>>>    \tau_n(X,x_0)       @>>>  \sigma_n(X,x_0,G)     @>>>   G    @>>> 1\\
@.        @VVV                 @VVV                   @|         \\
1 @>>>    \tau_{n-1}(X,x_0)   @>>>  \sigma_{n-1}(X,x_0,G) @>>>   G    @>>> 1,
\end{CD}
\end{equation}
where the first two vertical homomorphisms have sections.
Since this diagram is commutative for all $k$, it follows that the canonical projections
$\sigma_n(X,x_0,G)\to \sigma_{n-1}(X,x_0,G)$ and $\tau_n(X,x_0)\to \tau_{n-1}(X,x_0)$ have isomorphic
kernels, which we call $K$. But, the split exact sequence of Fox \cite{fox} takes the following form
(see \cite{ggw})
\begin{equation}\label{fox-exact}
0\to \tau_{n-1}(\Omega X,\ovl{x_0}) \to \tau_n(X,x_0) \stackrel{\dashleftarrow}{\to} \tau_{n-1}(X,x_0) \to 1.
\end{equation}
Thus,
$$K\cong \tau_{n-1}(\Omega X,\overline{x_0}),$$
where $\overline {x_0}$ denotes the constant loop at $x_0$. Now, consider the following commutative diagram
{\tiny
\begin{equation*}
\begin{CD}
1 @>>>    K   @>>>  \sigma_n(X,x_0,G) @>>>   \sigma_{n-1}(X,x_0,G) @>>> 1\\
@.    @VV{p_*}V           @VV{p_\ast}V                 @VV{\cong}V  \\
1 @>>>    \tau_{n-1}(\Omega (X/G),\overline {p(x_0)})   @>>>  \tau_n(X/G,p(x_0)) @>>>   \tau_{n-1}(X/G,p(x_0)) @>>> 1.
\end{CD}
\end{equation*}
}
Note that $\tau_*(\Omega Z)$ is a direct product of higher homotopy groups of $Z$. Since $G$ acts freely on
$X$, the orbit map $p:X\to X/G$ is a finite cover. Hence, the map $p_\ast : K\cong\tau_{n-1}(\Omega X,x_0)\to \tau_{n-1}(\Omega(X/G),p(x_0))$
is an isomorphism. Thus, by the Five Lemma and inductive hypothesis, the homomorphism $p_* : \sigma_n(X,x_0,G)\to \tau_n(X/G,p(x_0))$ is an isomorphism.
\end{proof}

Since, in \eqref{fox-exact}, the kernel of $\tau_n \to \tau_{n-1}$ is $\tau_{n-1}(\Omega X,\ovl{x_0})$, we have the following result generalizing \eqref{fox-exact}.

\begin{theorem}
The following sequence
\begin{equation}\label{rhodes-split}
0\to \tau_{n-1}(\Omega X,\ovl{x_0}) \to \sigma_n(X,x_0,G) \stackrel{\dashleftarrow}{\to} \sigma_{n-1}(X,x_0,G) \to 1
\end{equation}
is split exact.
\end{theorem}

\begin{remark}\label{sigma1-action}
Let $s_{n-1}:\sigma_{n-1} \to \sigma_n$ be the section. Consider the injective homomorphism $s_{n-1}\circ\cdots\circ s_2\circ s_1:\sigma_1 \to \sigma_n$ and the action $\Theta_n:\sigma_n \to {\rm Aut}(\sigma_n)$ given by conjugation. The action of $\sigma_1$ on $\sigma_n$ is given by
$$
\Theta_n\circ s_{n-1}\circ\cdots\circ s_1:\sigma_1\to {\rm Aut}(\sigma_n).
$$
Likewise, for the Fox torus homotopy groups, there is an action of $\pi_1=\tau_1$ on $\tau_n$. In fact, the action of
$\tau_1$ on $\tau_n$ is simply the restriction of that of $\sigma_1$ on $\sigma_n$. Furthermore, since $\tau_{n-1}(\Omega X,\ovl{x_0})\cong \prod^n \pi_i(X,x_0)^{\alpha_i}$, one can show that $\pi_n(X,x_0)$ is a (embedded) normal subgroup of $\sigma_n$. By embedding $\sigma_1$ in $\sigma_n$, one obtains an action of $\sigma_1(X,x_0,G)$ on $\pi_n(X,x_0)$.
\end{remark}

\section{Evaluation subgroups of Rhodes groups}

Given a space $X$, we defined in \cite{ggw} the Gottlieb-Fox groups to be the eva\-lu\-ation subgroups
$G\tau_n:=G\tau_n(X,x_0):={\rm Im} (ev_*:\tau_n(X^X,1_X)\to \tau_n(X,x_0))$
of the torus homotopy groups $\tau_n$ for $n\ge 1$, where $X^X:=\mbox{Map}(X,X)$
is the map space.
\par In this section, we define and study the analogous evaluation subgroups of the Rhodes groups $\sigma_n$ for $n\ge 1$.
\begin{definition}
Given a $G$-space $X$, consider the pointwise action of $G$ on the space $X^X$,
i.e., $(gf)(x):=gf(x)$ for $g\in G$, $f\in X^X$ and $x\in X$.
\par The evaluation subgroup $$G\sigma_n:=G\sigma_n(X,x_0,G):={\rm Im} (ev_*:\sigma_n(X^X,1_X,G)\to \sigma_n(X,x_0,G))$$
of $\sigma_n$ is called the $n$-th {\em Gottlieb-Rhodes group} of a $G$-space $X$.
\end{definition}
To relate the Gottlieb-Rhodes groups with the Gottlieb-Fox groups,
we consider the homomorphism $p_n:G\sigma_n \to G$ given by $[f;g]\mapsto g$ for $[f;g]\in G\sigma_n$.
\begin{theorem}\label{gottlieb-rhodes-gottlieb-fox}
The following sequence
\begin{equation} \label{evaluation-exact}
1 \to G\tau_n \to G\sigma_n \stackrel{p_n}{\to} G_0 \to 1
\end{equation}
is exact. Here, $G_0$ is the subgroup of $G$ consisting of elements $g$ considered as homeomorphisms of $X$
which are freely homotopic to $1_X$.
\end{theorem}
\begin{proof}
Since $[f;g]\in G\sigma_n$, we get $[f;g]=ev_*([F;g])$ for some $[F;g]\in \sigma_n(X^X,1_X,G)$.
This means that $g_*\sim 1_X$ (or equivalently $g\in G_0$) via the homotopy determined by the restriction map
$F_| :I\times\{1,\ldots,1\}\to X^X$. Thus, we get $p_n(G\sigma_n)=G_0$.
By the naturality of the evaluation map and the fact that $\mbox{Ker}(\,p_n\subseteq \tau_n$,
it follows that $\mbox{Ker}\,p_n=G\tau_n$.
\end{proof}

The group $G\sigma_n$ was already defined by M.\ Woo and Y.\ Yoon \cite{wy1} who asked whether $G\sigma_n$ is abelian in general. In view of \eqref{evaluation-exact}, we expect that $G\sigma_n$ is non-abelian in general.
\begin{example}\label{ex1}
Let $X=\mathbb RP^3=\mathbb{S}^3/\{\pm I\}$ be the $3$-dimensional real projective space and
$G=\mathbb Z_2 \oplus \mathbb Z_2$. Consider the  free action of the quaternionic group $Q_8 \subset \mathbb{S}^3$
on the sphere $\mathbb{S}^3$. Then, $G$ acts on $X$ as a quotient of $Q_8$. Then, $G\sigma_1=Q_8$, $G\tau_1=\mathbb Z_2$ and $G_0=G=\mathbb Z_2 \oplus \mathbb Z_2$.
\end{example}

We have three evaluation subgroups, namely the classical Gottlieb groups
$G_n(X)=G_n(X,x_0)$, the Gottlieb-Fox groups $G\tau_n(X)=G\tau_n(X,x_0)$
and the Gottlieb-Rhodes groups $G\sigma_n(X,G)=G\sigma_n(X,x_0,G)$, where
$(X,G)$ is a transformation group. We shall compare these different
notions. Recall that a space $X$ is a {\em Gottlieb space} if
$G_n(X)=\pi_n(X)=\pi_n(X,x_0)$ for all $n\ge 1$. Similarly, we say that
$X$ is {\it Gottlieb-Fox} or {\it Gottlieb-Rhodes} for a $G$-space $X$ if
$G\tau_n(X)=\tau_n(X)$ and $G\sigma_n(X,G)=\sigma_n(X,G)$, respectively
for all $n\ge 1$. Certainly, any $H$-space is Gottlieb-Fox.

\begin{remark} \label{rmk} It is straightforward to see, using Theorem \ref{gottlieb-rhodes-gottlieb-fox},
that a $G$-space $X$ is Gottlieb-Rhodes if and
only if $X$ is Gottlieb-Fox and $G_0=G$.
\end{remark}

\begin{proposition}\label{gottliebfox=fox}
{\em Suppose that $X$ is a compactly generated Hausdorff path connected space.
Then, $X$ is Gottlieb if and only if it is a Gottlieb-Fox space.}
\end{proposition}
\begin{proof} It was shown in \cite{ggw} that the Gottlieb-Fox groups are direct products of the classical
Gottlieb groups. More precisely, there is an isomorphism
$$
G\tau_n(X)\cong \prod_{i=1}^n G_i(X)^{\gamma_i},$$
where $\gamma_i=\binom{n-1}{i-1}$. From \cite{fox}, we also have
$$
\tau_n(X)\cong \prod_{i=2}^n \pi_i(X)^{\alpha_i} \rtimes \tau_{n-1}(X),
$$
where $\alpha_i=\binom{n-2}{i-2}$. Suppose $X$ is Gottlieb. Since $\tau_1=\pi_1$ and $G\tau_1=G_1$, we have $\tau_1=G\tau_1$. By inductive hypothesis, we may assume that $G\tau_n(X)=\tau_n(X)$. Now,
\begin{equation*}
\begin{aligned}
\tau_{n+1}(X)&\cong \prod_{i=2}^{n+1} \pi_i(X)^{\beta_i} \rtimes G\tau_n(X) \\
             &\cong \prod_{i=2}^{n+1} \pi_i(X)^{\beta_i} \rtimes \prod_{i=1}^n \pi_i(X)^{\gamma_i},
\end{aligned}
\end{equation*}
where $\beta_i=\binom{n-1}{i-2}$. Since the semi-direct is given by the Whitehead product and $X$ is Gottlieb, the action is trivial and thus
$$
\tau_{n+1}(X) \cong \prod_{i=2}^{n+1} \pi_i(X)^{\beta_i} \times \prod_{i=1}^n \pi_i(X)^{\gamma_i}.$$ For $2\le i\le n$, $$\beta_i+\gamma_i=\binom{n-1}{i-2} + \binom{n-1}{i-1} = \binom{n}{i-1}.$$
It is easy to see that
$$
\prod_{i=1}^{n+1} G_i(X)^{\binom{n}{i-1}}\cong G\tau_{n+1}(X)=\tau_{n+1}(X) \cong \prod_{i=1}^{n+1} \pi_i(X)^{\binom{n}{i-1}}.
$$
Thus, by induction, $X$ is a Gottlieb-Fox space.

Conversely, if $X$ is Gottlieb-Fox, then $G\tau_1(X)=\tau_1(X)$ is equivalent to $G_1(X)=\pi_1(X)$. Suppose
$G_i(X)=\pi_i(X)$ for $i\le n$. Then, $\tau_{n+1}(X)\cong \prod_{i=1}^{n+1}\pi_i(X)^{\beta_i} \times \tau_n(X)$ is
a direct product since $\tau_n(X)=G\tau_n(X)\cong \prod_{i=1}^n G_i(X)^{\gamma_i}$ and Whitehead products
with Gottlieb elements vanish. Now,
$$
\prod_{i=1}^n G_i(X)^{\binom{n}{i-1}}
\times G_{n+1}(X) \cong G\tau_{n+1}(X)=\tau_{n+1}(X) \cong \prod_{i=1}^n \pi_i(X)^{\binom{n}{i-1}} \times \pi_{n+1}(X).$$
But, by inductive hypothesis, $\pi_i(X)=G_i(X)$ for $i\le n$. Note that the equality above mens that the inclusion $G\tau_{n+1}(X)\hookrightarrow \tau_{n+1}(X)$ is an isomorphism. This inclusion is given by the product of the respective inclusions $G_{n+1}(X) \subseteq \pi_{n+1}(X)$ and $G_i(X)\subseteq \pi_i(X)$ for $i\le n$. Since the inclusions $G_i(X) \subseteq \pi_i(X)$, for $i\le n$, are isomorphisms, it follows that the inclusion $G_{n+1}(X)\hookrightarrow \pi_{n+1}(X)$ must be an isomorphism. Hence, $X$ is a Gottlieb space.
\end{proof}

\begin{example}\label{ex2}
Let $X=\mathbb{S}^n$ be the $n$-sphere and $G$ be a finite group acting freely on $X$. If $n$ is even, then
$G=\mathbb Z_2$ in which case $G\ne G_0=\{e\}$. If $n$ is odd, then for any $g\in G$, the Lefschetz number
$L(g_*)$ is zero since $g$ has no fixed points. It follows that $L(g_*)=1-\deg g_*=0$ and this
implies that $\deg g_\ast=1$. Thus, $g_\ast$ is homotopic to $1_X$. In this case, $G=G_0$.
\end{example}

\begin{example}\label{ex3}
Let $X=T^3$ be the 3-dimensional torus and $G=\mathbb Z_2=\langle t\rangle$. The action of $G$ on $X$
is given by $t\cdot (a,b,c) \mapsto (Aa,\bar b, \bar c)$, where $A$ is the antipodal map on $\mathbb{S}^1$
and $\bar z$ denotes the complex conjugate of $z$. The map induced by the generator $t$ is of degree $1$ but
is not homotopic to the identity so that $G\ne G_0$. It is well known that $H$-spaces are Gottlieb spaces.
This example shows that $X$ being Gottlieb-Fox does not imply $X$ that is Gottlieb-Rhodes without the condition
$G=G_0$.
\end{example}

Next, we generalize a result of Gottlieb \cite{gottlieb} on the Gottlieb subgroups in a fibration. Let $G$ be a finite group. By a $G$-fibration, we mean a $G$-map $p:E\to B$ satisfying the $G$-Covering Homotopy Property for any $G$-spaces.

\begin{theorem}\label{G-fibration}
Let $G$ be a finite group and $E \stackrel{p}{\to} B$ be a $G$-fibration such that $B^G\ne \emptyset$. Let $F=p^{-1}(b_0)$ for some $b_0\in B^G$ and choose $\tilde b_0\in F$. Then, there exists a homomorphism $d_*:\sigma_n(\Omega B, \ovl{b_0},G)\to \sigma_n(F,\tilde b_0,G)$ such that
$$
d_*(\sigma_n(\Omega B,\ovl{b_0},G)) \subseteq G\sigma_n(F,\tilde b_0,G).
$$
\end{theorem}
\begin{proof} Let $\Omega B$ be the space of loops based at $b_0$. Consider the $G$-action on
$\Omega B\times G$ given by $g\cdot (\lambda, g')=(g\lambda, gg')$ for $(\lambda,g')\in\Omega B\times G$ and $g\in G$.
Because the definition of $\sigma_n$ automatically restricts to the component of
the basepoint and the group $G$ is finite, the projection $p_G:\Omega B\times G \to \Omega B$
induces an isomorphism
$$
{p_G}_*:\sigma_n(\Omega B\times G, (\ovl{b_0},e),G) \to \sigma_n(\Omega B,\ovl{b_0},G),
$$
where $e\in G$ denotes the identity element.

Using the $G$-Covering Homotopy Property, there exists a $G$-lifting function $\mu_G:\Omega B\times G \to F^F$.
Define a map $d_*$ to be the composite
$$
\sigma_n(\Omega B,\ovl{b_0},G) \stackrel{{p_G}_*^{-1}}{\to} \sigma_n(\Omega B\times G, (\ovl{b_0},e),G) \stackrel{{\mu_G}_*}{\to} \sigma_n(F^F,1_F,G) \stackrel{{ev}_*}{\to} \sigma_n(F,\tilde{b_0},G).
$$
Then, the assertion follows from the definitions of $d_*$ and of $G\sigma_n$.
\end{proof}

\begin{remark}
When $G=\{1\}$, the map $d_*$ becomes the boundary homomorphism induced by the action of $\Omega B$ on $F$
so that Theorem \ref{G-fibration} becomes the first part of \cite[Theorem 2.2]{ggw}.
\end{remark}

\section{Equivariant Gottlieb groups and orbit spaces}

In \cite{gg}, equivariant Gottlieb groups $\{{\bf G}_n(X^H,x_0)\}$ or simply $\{{\bf G}_n(X^H)\}$ were defined
for every subgroup $H\le G$ of a $G$-space $X$, where $X^H$ is the subspace of $X$ given by $H$-invariant
points. These groups are subgroups of the classical homotopy groups
$\pi_n(X^H)$. When $n=1$, $\{{\bf G}_1(X^H)\}$ are the same as the $G$-Jiang subgroups as defined in
\cite{wong}. In this section, we investigate relationships among the various evaluation groups.

For any positive integer $n$, a space $X$ is {\it $n$-Gottlieb} if $G_n(X)=\pi_n(X)$. Similarly, $X$ is {\it $n$-Gottlieb-Fox} if $G\tau_n(X)=\tau_n(X)$. Recall from \cite{gg} that for any $n\ge 1$ and for any subgroup $H\le G$,
$${\bf G}_n(X^H)=
{\rm Im}(ev_*:\pi_n({\rm Map}_{WH}(X^H,X^H),1_{X^H}) \to \pi_n(X^H,x_0)),$$
where ${\rm Map}_{WH}(X^H,X^H)$
is the space of $WH=NH/H$-maps on $X^H$, where $NH$ is the normalizer of $H$ in the group $G$ and $WH$ is the Weyl group. For any positive integer $n\ge 1$, a $G$-space $X$ is said to be
{\em $n$-Gottlieb-Rhodes} if $G\sigma_n(X,x_0,G)=\sigma_n(X,x_0,G)$. Similarly, we say that $X$ is
{\em equivariant $n$-Gottlieb} if ${\bf G}_n(X^H)=\pi_n(X^H)$ for every subgroup $H$. While one can define
analogously the notion of {\it equivariant $n$-Gottlieb-Rhodes}, i.e., $G\sigma_n(X^H,x_0,WH)=\sigma_n(X^H,x_0,WH)$ for every subgroup $H$, for the rest of the paper we will only study the various evaluation subgroups when the $G$-action is free.

\begin{theorem}\label{orbit-space}
Suppose a group $G$ acts freely on a space $X$ with a basepoint $x_0$.
\begin{enumerate}
 \item For $n\ge 2$, if $X$ is equivariant $n$-Gottlieb then $X/G$ is $n$-Gottlieb.
 \item For $n\ge 1$, if $X/G$ is $n$-Gottlieb then $X$ is equivariant $n$-Gottlieb.
 \item Suppose $X$ is a finite aspherical $G$-CW space. If $X$ is equivariant $1$-Gottlieb then $X/G$ is $1$-Gottlieb.
\end{enumerate}
\end{theorem}
\begin{proof} Certainly, given a free $G$-action on $X$, the subspaces $X^H$ are empty for any  non-trivial
subgroup $H$ of $G$ and $X^{\{e\}}=X$.

(1) For the obvious map $q : \mbox{Map}_G(X,X)\to \mbox{Map}(X/G,X/G)$, consider the following commutative diagram
\begin{equation*}
\begin{CD}
\pi_n({\rm Map}_{G}(X,X),1_X)   @>{ev_*}>>  \pi_n(X,x_0)\\
@V{q_*}VV                         @V{p_*}VV  \\
\pi_n({\rm Map}(X/G,X/G),1_{X/G})   @>{ev_*}>>  \pi_n(X/G,p(x_0))
\end{CD}
\end{equation*}
for any $n\ge 2$. If $X$ is equivariant $n$-Gottlieb then the evaluation
map $ev_*:\pi_n({\rm Map}_G(X,X),1_X)\to \pi_n(X,x_0)$ must be onto.
Furthermore, since $p_*:\pi_n(X,x_0)\to \pi_n(X/G,p(x_0))$ is an
isomorphism, the commutativity of the diagram implies that
$ev_*: \pi_n({\rm Map}(X/G,X/G),1_{X/G}) \to \pi_n(X/G,p(x_0))$ is also
surjective, i.e., $G_n(X/G;p(x_0))=\pi_n(X/G,p(x_0))$.

(2) The case when $n=1$ was already proven in \cite[Proposition 4.9]{wong}.
For $n\ge 2$, it was proven in \cite[Proposition 3.3]{gg} that $p_*({\bf G}_n(X))\cong G_n(X/G)$. Since $p_*:\pi_n(X)\to \pi_n(X/G)$ is an isomorphism for $n\ge 2$ and $G_n(X/G)=\pi_n(X/G)$ by assumption, it follows that ${\bf G}_n(X)=\pi_n(X)$, i.e., $X$ is equivariant $n$-Gottlieb.

(3) Since $X$ is a finite aspherical $CW$ space, so is $X/G$. The fundamental groups $\pi_1(X)$ and
$\pi_1(X/G)$ are finitely generated and torsion-free. Since $X$ is equivariant $1$-Gottlieb, $\mbox{\bf G}_1(X)=\pi_1(X)$
and consequently $\pi_1(X)\cong \mathbb Z^d$ for some integer $d>0$. By \cite[Proposition 3.3]{gg},
$p_*({\bf G}_1(X))\subseteq G_1(X/G)$. Moreover, a result of Gottlieb asserts that $G_1(X/G)$ is the center
of $\pi_1(X/G)$. It then follows that $p_\ast(\pi_1(X))=p_*({\bf G}_1(X)) \subseteq G_1(X/G)$ is central in $\pi_1(X/G)$.
In other words, the exact sequence
$$
0\to \mathbb Z^d \cong \pi_1(X) \to \pi_1(X/G) \to G \to 1
$$
determined by the orbit map $p : X\to X/G$ is a central extension. By \cite{vasquez}, the group $\pi_1(X/G)\cong \mathbb Z^d$ is free abelian
of the same rank as $\pi_1(X)$. Then, since $X/G$ is aspherical, $G_1(X/G)=\mbox{Center}(\pi_1(X/G)=\pi_1(X/G))$
and the space $X/G$ is $1$-Gottlieb.
\end{proof}

For $n=1$, (1) of Theorem \ref{orbit-space} does not hold as we show in the following examples.

\begin{example}\label{ex4}
Let $X=\mathbb{S}^3\times \mathbb{S}^3 \times \mathbb{S}^3$, $G=\mathbb Z_2=\langle t\rangle$, where the free
action of $G$ on $X$ is given by $t\cdot (x,y,z)=(-x,\ovl y, \ovl z)$. Here, $\ovl w$ denotes the conjugate
of $w$ in $\mathbb{S}^3$ which is regarded as the unit quaternions. Note that ${\bf G}_1(X)=\{1\}$ since
$\pi_1(X)=\{1\}$, and $G_1(X/G)=\{1\}$ while $\pi_1(X/G)=\mathbb Z_2$. To see this, observe that the elements
of $G_1(X/G)$ act trivially on $\pi_3(X/G)\cong\pi_3(X)$ but the automorphism $t_{\sharp} \in {\rm Aut}(\pi_3(X))$
for $t_{\sharp}(1\oplus 1\oplus 1)=(1\oplus (-1)\oplus (-1))$ is non-trivial.
\end{example}

\begin{example}\label{ex5}
Let $X=\mathbb{S}^3$ and $G$ be a finite subgroup of $\mathbb{S}^3$. The free action of $G$ on $X$ is multiplication
in $\mathbb{S}^3$ so that $X/G$ is the coset space $\mathbb{S}^3/G$. By \cite{lang},
$G_n(\mathbb{S}^3/G)\cong G_n(\mathbb{S}^3)$ for all $n\ge 1$. Since $\mathbb{S}^3$ is a topological group,
it is $n$-Gottlieb for all $n\ge 1$, i.e., $\mathbb{S}^3$ is a Gottlieb space. Using
\cite[Proposition 3.3]{gg}, we have ${\bf G}_n(\mathbb{S}^3)\xrightarrow[\cong]{{p_*}} G_n(\mathbb{S}^3/G)$
for $n\ge 2$ so that ${\bf G}_n(\mathbb{S}^3)\cong G_n(\mathbb{S}^3)\cong \pi_n(\mathbb{S}^3)$. For $n=1$,
${\bf G}_1(\mathbb{S}^3)=\{1\}=\pi_1(\mathbb{S}^3)$. Therefore, $\mathbb{S}^3$ is
equivariant $n$-Gottlieb for all $n\ge 1$. However, if $G$ is not abelian, then $\mathbb{S}^3/G$ is not $1$-Gottlieb.
\end{example}

\begin{remark} J.\ Oprea \cite{oprea} proved that if $G$ is a finite group acting freely on $\mathbb{S}^{2n+1}$ then $G_1(\mathbb{S}^{2n+1}/G)$
is the center of $G$. In particular, when $n=1$, Oprea's result asserts that $\mathbb{S}^3/G$ is Gottlieb if
and only if $G$ is abelian. This result cannot be generalized to arbitrary simply connected Gottlieb spaces. In fact, Example \ref{ex4} serves as a counter-example in which $X$ is a simply-connected Gottlieb space (since it is a topological group) admitting a free $G=\mathbb Z_2$ action but $X/G$ is not a Gottlieb space.
\end{remark}

J.\ Siegel \cite{siegel} gave the first example of a finite dimensional Gottlieb space that is not an $H$-space. Following our discussion in this section, one can construct similar but even simpler examples as follows.

\begin{example}\label{ex6}
Let $X=\mathbb{S}^3/\mathbb Z_4$ be the coset space of $\mathbb{S}^3$ by the finite subgroup $\mathbb Z_4$. Then, $X$ is a Gottlieb space but not an $H$-space. The fact that $X$ is a Gottlieb space follows from the same argument used in Example \ref{ex5} except now $G=\mathbb Z_4$ is abelian so that $X$ is indeed $1$-Gottlieb.
The fact that $X$ does not admit an $H$-structure follows from \cite[Theorem 4, p.\ 87]{borel}.
\end{example}

We end this section by studying the Gottlieb-Rhodes groups of a free $G$-space and the Gottlieb groups of the orbit space when the space is aspherical.

\begin{proposition}\label{G-gottlieb-gottlieb-rhodes}
{\em Suppose a finite group $G$ acts freely on a finite aspherical $G$-$CW$ space $X$
with a basepoint. Then $X/G$ is Gottlieb if and only if $X$ is a Gottlieb-Rhodes space.}
\end{proposition}
\begin{proof}
By Theorem \ref{orbit-space}, $X/G$ is Gottlieb if and only if
$X$ is equivariant Gottlieb. Since the $G$-action is free, an equivariant Gottlieb space is the same as
a Gottlieb $G$-space. Using Remark \ref{rmk} and Proposition \ref{gottliebfox=fox}, $X$ is Gottlieb-Rhodes if and only if $X$ is Gottlieb and
$G_0=G$. Thus, if $X$ is Gottlieb-Rhodes then $X/G$ is Gottlieb.
\par To prove the converse, it suffices to show
that if $X/G$ is Gottlieb then $G_0=G$. Since $X/G$ is Gottlieb the groups $\pi_1(X/G)$ and $\pi_1(X)$ are
abelian. Hence, because the space $X$ is path connected, there is a canonical induced $G$-action on $\pi_1(X)$.
But $p:X \to X/G$ is a finite cover, so it follows that the short exact sequence $0\to \pi_1(X) \to \pi_1(X/G) \to G \to 1$
is a central extension so that $G$ acts trivially on $\pi_1(X)$.
That is, every $g\in G$ induces the identity homomorphism on $\pi_1(X)$. Since $X$ is aspherical, it follows that $g$ is homotopic to the identity map
$1_X$. This means that $G_0=G$.
\end{proof}

\begin{remark} In general, $X$ being equivariant Gottlieb does not
imply that $X$ is Gottlieb-Rhodes even when $G$ acts freely on $X$.
In fact, Example \ref{ex4} furnishes such an example since
$G_0=\{1\}\ne G=\mathbb Z_2$ and $X$ is equivariant Gottlieb.
\end{remark}


\end{document}